\documentstyle{article}
\newtheorem{theorem}{Theorem}[section]

\newtheorem{Prop}{Proposition}[section]
\newcommand{\si}{\sigma}
\newcommand{\Si}{\Sigma}
\newcommand{\mbf}{\mathbf}
\newcommand{\dsty}{\displaystyle}
\newcommand{\dsum}{\displaystyle\sum}
\newcommand{\dlim}{\displaystyle\lim}
\newcommand{\ra}{\rightarrow}
\newcommand{\Ra}{\Rightarrow}

\newcommand{\al}{\alpha}
\newcommand{\be}{\beta}
\newcommand{\de}{\delta}
\newcommand{\ep}{\epsilon}

\newcommand{\tri}{\triangle}

\begin{document}
\title{Exponential Random Energy Model}
\author{Nabin Kumar Jana\footnote{Stat-Math Unit, Indian Statistical Institute, 203 B. T. Road, Kolkata,
India.
e-mail: nabin\_r@isical.ac.in}\\
 {\small Indian Statistical Institute, Kolkata, India}}

\date{March 4, 2005 (Revised May 1, 2005)}
\pagestyle{myheadings} \markboth{Jana N. K.}{Exponential Random
Energy Model}

\maketitle
\begin{abstract}
In this paper the Random Energy Model(REM) under exponential type
environment is considered which includes double exponential and Gaussian cases.
 Limiting Free Energy is evaluated in these models. Limiting Gibbs' distribution
 is evaluated in the double exponential case.\\

\end{abstract}
\date{}
\maketitle

\section{Introduction}

Usually the Random Energy Model$^{\cite{D}}$ is considered in a
Gaussian$^{\cite{B,CEGG,DW,OP,T1}}$ environment. In this paper we
discuss the same under a double exponential environment. It is
interesting to note that in our analysis the distribution of
Hamiltonian, $H_N$ does not depend on $N$. We use large deviation
method$^{\cite{DW}}$ to calculate the limiting free energy. There is
a phase transition at $\be=1$. The methods carry over to a more
general exponential type family that includes the Gaussian case as
well and we provide explicit formulae for the free energy.

We use Talagrand's$^{\cite{T1}}$ approach to obtain the limiting
Gibbs' distribution in the low temperature regime. It is interesting
to note that the limit is again a Poisson-Dirichlet distribution.
Observe that when $X$ is double exponential with parameter one,
${\mbf{E}}(e^{\be X})$ does not exist for $\be >1$. For $0<\be < 1$
we obtain -- as expected -- uniform distribution as the infinite
volume limit of Gibbs' distributions. This is done via an
interesting variant of the strong law of large numbers. These
methods carry over to the Gaussian case as well.

\section{Free energy}
For each configuration $\si \in \Si_N = \{-1, 1\}^N$ of an $N$
particle system, the Hamiltonian is $H_N(\si)$. Since we are
considering REM, $\{H_N(\si)\}$ are i.i.d. In fact, we assume that
they are double exponential, that is, have density (not depending on
$N$)
$$\phi_N(x) = \frac{1}{2} e^{-|x|}, \vspace{3ex} -\infty < x < \infty.$$

The partition function of the system is $$Z_N(\be) = \dsum_\si
e^{-\be H_N(\si)}=2^N{\mbf{E}}_\si e^{-N\be \frac{H_N(\si)}{N}},$$
where $\be >0$ is the inverse temperature and $\mbf{E_\si}$ stands
for expectation w.r.t $\si$ when $\Si_N$ has uniform distribution.
Hence the free energy of the system is \[\frac{1}{N} \log Z_N(\be)=
\log 2 + \frac{1}{N} \log{\mbf{E}}_{\si} e^{-N\be
\frac{H_N(\si)}{N}}.\] Now let $\mu_N$ be the induced (random)
probability on ${\mbf{R}}$ via the map
\[ \si \mapsto
\frac{H_N(\si)}{N}\] when $\Si_N$ has uniform distribution.

\begin{Prop}
$\mu_N \Ra \de_0$ a.s. as $N\ra \infty$.

That is, for almost every sample point, the sequence of random
measures $\{\mu_N\}$ converges weakly to point mass at $0$.
\end{Prop}


\noindent {\it Proof:} For any $\ep>0$, define $\tri(\ep) =
[-\ep,\ep] \subset {\mbf{R}}.$ Now by Markov inequality,
$${\mbf{P}}(\mu_N(\tri^c(\ep))>\ep) < \frac{1}{\ep}{\mbf{E}}\mu_N(\tri^c(\ep) )< \frac{1}{\ep}{\mbf{P}}(|H_N|>\ep N) = \frac{1}{\ep}e^{-\ep N}.$$

Apply Borel-Cantelli.\\

Let $\tri \subseteq {\mbf{R}}$ be an open interval. Put $m =
\dsty\inf_{x \in \tri} |x|$ and $M = \dsty\sup_{x \in \tri}|x|$, and
$q_N = P(\frac{H_N}{N} \in \tri)$. These quantities, of course,
depend on $\tri$. Observe that $\tri \subseteq (-M, -m] \cup [m,
M)$, so that
\begin{equation}
q_N \leq \dsty\int_{Nm}^{NM} e^{-x} dx \leq \dsty\int_{Nm}^{\infty}
e^{-x} dx = e^{-Nm}
\end{equation}
and
\begin{equation}
q_N \geq \frac{1}{2} \dsty\int_{Nm}^{NM} e^{-x} dx > \frac{1}{2}
\dsty\int_{Nm}^{Nm + \de} e^{-x} dx > \frac{\de}{2} e^{-(Nm + \de)},
\end{equation}
for any $\de$, $(0<\de < M - m)$. Both (1) and (2) remain true even
if $m=0.$

\begin{Prop}
If $m > \log2$, then almost surely eventually $\mu_N (\tri) = 0$.

Hence almost surely eventually, $\mu_N[-\log2,\log2]=1$
\end{Prop}

\noindent {\it Proof:} By definition,
$$\mu_N(\tri) =
\frac{1}{2^N}\dsty\sum_{\si}{\mbf{1}}_{\frac{H_N(\si)}{N}}(\tri)$$

Hence $\{\mu_N(\tri) >0\} =
\{\dsty\sum_{\si}{\mbf{1}}_{\frac{H_N(\si)}{N}}(\tri) > 1\}$. Now by
Chebyscheff's inequality, $ {\mbf{P}}(\mu_N(\tri) >0) <
{\mbf{E}}\dsty\sum_{\si}{\mbf{1}}_{\frac{H_N(\si)}{N}}(\tri) =2^N
q_N.$

Since $m> \log2$, (1) implies, $\dsty\sum_{N\geq 1}
{\mbf{P}}(\mu_N(\tri) >0) < \infty.$ Borel-Cantelli completes the
proof.

\begin{Prop}
If $m < \log2$, then for any $\ep
> 0$ a.s. eventually
\[(1-\ep)q_N \leq \mu_N(\tri) \leq (1+\ep)q_N.\]
\end{Prop}
\noindent {\it Proof:} Note that
\[\begin{array}{lll}
var(\mu_N(\tri))&=& {\mbf{E}}(\mu_N(\tri))^2 - [{\mbf{E}}\mu_N(\tri)]^2\\
&=& \frac{1}{2^{2N}}\dsty\sum_{\si, \tau} \left[{\mbf{E}} {\mbf{1}}_{\frac{H_N(\si)}{N}}(\tri){\mbf{1}}_{\frac{H_N(\tau)}{N}}(\tri)- q_N^2\right]\vspace{1ex}\\
&=& \frac{1}{2^{2N}}\dsty\sum_{\si= \tau} \left[{\mbf{E}} {\mbf{1}}_{\frac{H_N(\si)}{N}}(\tri){\mbf{1}}_{\frac{H_N(\tau)}{N}}(\tri)- q_N^2\right]\vspace{1ex}\\
&\leq& \frac{1}{2^{2N}}\dsty\sum_{\si} {\mbf{E}}{\mbf{1}}_{\frac{H_N(\si)}{N}}(\tri)\\
&=& \frac{q_N}{2^N}
\end{array}\]
Hence for any $\ep >0$, by Chebycheff's inequality
$$
{\mbf{P}}\left[|\mu_N(\tri)-{\mbf{E}}\mu_N(\tri)|>\ep
{\mbf{E}}\mu_N(\tri)\right] < \frac{1}{\ep^2 2^{N}q_N}.$$

Using (2) and the fact $m <\log2$, we get $\sum_N 2^{-N}q_N^{-1}<
\infty$, so that $\sum_N{\mbf{P}}(|\mu_N(\tri)-q_N|>\ep q_N)
<\infty.$ Borel-Cantelli
completes the proof.\\

Propositions 2.3 and 2.2 combined with the inequalities (1) and (2)
yield,
\begin{Prop}
Almost surely,
\[\begin{array}{llll}
\dsty\lim_{N\ra \infty} \frac{1}{N} \log\mu_N(\tri) &=& - m & \mbox{
if } m < \log2\\
&=& -\infty & \mbox{ if } m > log2.
\end{array}\]
\end{Prop}

Now let us consider the map $I: {\mbf{R}} \ra {\mbf{R}}^+$, defined
as follows,
\[\begin{array}{llll}
I(x) &=&  |x| &\mbox{ if } -\log2 \leq x \leq \log2\\
&=& \infty &\mbox{ otherwise}.
\end{array}\]
\begin{theorem}
Almost surely, the sequence $\{\mu_N\}$ satisfies the large
deviation principle with rate function $I$.
\end{theorem}
\noindent {\it Proof:} The collection $\mathcal{A}$ of open
intervals with rational end points is a countable base for
${\mbf{R}}$. For $\tri \in \mathcal{A}$, put $L_\tri = -
\dsty\lim_{N\ra \infty} \frac{1}{N} \log \mu_N(\tri)$. Note that,
for $x\in{\mbf{R}}$ Proposition 2.4 implies, $I(x) = \dsty\sup_{x
\in \tri \in \mathcal{A}} L_\tri.$

Since almost surely $\{\mu_N\}$ is supported on a compact set,
Theorem 4.1.11 of Dembo and Zeitouni$^{\cite{DZ}}$ completes the
proof.

\begin{theorem}
For almost every sample point,
\[\begin{array}{llll}
\dsty\lim_{N \ra \infty} \frac{1}{N} \log Z_N(\be) &=& \log2
&\mbox{if } \be \leq 1\\
&=&\be \log2 &\mbox{if }  \be \geq 1.\\
\end{array}\]
\end{theorem}
\noindent {\it Proof:} By Theorem 2.1, almost surely, the sequence
$\{\mu_N\}$ satisfies large deviation principle with rate function
$I$. By Proposition 2.2, the sequence $\{\mu_N\}$ is supported on a
compact set. Varadhan's lemma with $h(x) = \be x$, $-\infty< x <
\infty$ gives
$$\begin{array}{lll}
\dsty\lim_{N\ra \infty}\frac{1}{N}\log Z_N(\be) &=&\log2 -
\dsty\inf_{x\in {\mbf{R}}}\{h(x)+I(x)\}\vspace{1ex}\\
&=&\log2 - \dsty\inf_{|x|\leq \log 2}\{\be x + |x|\}\vspace{1ex}\\
&=&\log2 - \dsty\inf_{0\leq x\leq \log 2}(1-\be) x.
\end{array}$$
Hence a.s. $$\begin{array}{llll} \dsty\lim_{N\ra
\infty}\frac{1}{N}\log Z_N(\be) &=& \log2 &\mbox{ if
} \be \leq 1\\
&=& \be \log2 &\mbox{ if } \be \geq 1
\end{array}$$
\vspace{1ex}

\noindent{\bf Remark 2.1} It is worth noting that the above
consideration hold for a general class of distributions. More
precisely, let $\al \geq 1$ be fixed. For each $N$, consider
$\{H_N(\si), \si \in \Si\}$ to be i.i.d. with density
$$\phi_N(x) = C_{\al,N}\; e^{-\frac{|x|^\al}{\al N^{\al-1}}}, \quad
-\infty <x<\infty,$$ where
$C_{\al,N}=\frac{1}{2\Gamma(\frac{1}{\al})}\left(\frac{\al}{
N}\right)^{\frac{\al-1}{\al}}$.

Of course, when $N=1$, this reduces to the case considered above and
for $N=2$ this becomes the Gaussian case usually considered in the
literature. For $\al>1$, similar calculations as above lead to the
rate function
$$\begin{array}{llll}
I(x) &=&  \frac{|x|^\al}{\al} &\mbox{ if } -(\al\log2)^{\frac{1}{\al}} \leq x \leq (\al\log2)^{\frac{1}{\al}}\\
&=& \infty &\mbox{ otherwise}.
\end{array}$$

and almost surely, the limiting free energy is
$$\begin{array}{llll}\dlim_{N\ra \infty} \frac{1}{N}\log Z_N(\be) &=&
\log2 + \frac{\al-1}{\al} \be^{\frac{\al}{\al-1}}& \mbox{ if } \be
\leq (\al \log2)^{\frac{\al-1}{\al}}\\
&=&\be(\al \log2)^{\frac{1}{\al}}& \mbox{ if } \be > (\al
\log2)^{\frac{\al-1}{\al}}
\end{array}.$$

For $\al =2$, this coincides with the known formula$^{\cite{T1}}$.
Of course, when $\al=1$ the formula, interpreted in the limiting
sense, is the one obtained earlier.

\section{{\bf Gibbs' Distribution}}
We return to the double exponential environment. Recall that Gibbs'
distribution for the $N$ particle system is the (random) probability
on $\Si_N$ defined as
$$G_N(\si) = \frac{e^{-\be H_N(\si)}}{Z_N(\be)},\quad \si \in
\Si_N.$$ We show that for $\be< 1$, the (random) Gibbs' distribution
$G_N$ converges weakly to the uniform probability on $\{-1,
1\}^{\infty}$ almost surely. Uniform probability here means the
product probability on $\{-1, 1\}^{\infty}$ where each coordinate
space has $(\frac{1}{2}, \frac{1}{2})$ probability. Since for each
$N$, $G_N$ is defined on $\{-1, 1\}^N$ the notion of convergence
here is to be carefully understood. This is made precise in Theorem
3.1 below.
\begin{theorem}
Fix $\be < 1$. Then almost surely, for any $K \geq 1$ and any $\si
\in \{-1, 1\}^K$, $\rho_N(\si) \ra \frac{1}{2^K}$ as $N \ra \infty$,
where $\rho_N$ is the marginal of $G_N$ on $\{-1, 1\}^K$.
\end{theorem}
\noindent {\it Proof:} For $0<\be<1$, define $Z_N'(\be)=\sum_\si
e^{-\be H_N(\si)}{\mbf{1}}_{\{-H_N(\si)\leq \de N\}}$, where $\de$
will be chosen latter depending on $\be$.

Since ${\mbf{P}}\{-H_N(\si)> \de N\} \leq \frac{1}{2}e^{-\de N}$,
for any $\de>log2$, Borel-Cantelli implies that almost surely,
eventually
\begin{equation}
Z_N(\be)=Z_N'(\be).
\end{equation}

Argument as in Proposition 2.3 and symmetry of distribution of $H_N$
lead to,
\begin{equation}{\mbf{P}}\left[|Z_N'(\be) -
{\mbf{E}}Z_N'(\be)|
> \ep {\mbf{E}}Z_N'(\be)\right] < \frac{{\mbf{E}}e^{2\be
H_N}{\mbf{1}}_{\{H_N\leq \de N\}}}{\ep^2 2^N ({\mbf{E}}e^{\be
H_N}{\mbf{1}}_{\{H_N\leq \de N\}})^2}.
\end{equation}

But, \begin{equation}{\mbf{E}}e^{\be H_N}{\mbf{1}}_{\{H_N\leq \de
N\}}>\frac{1}{1+\be}.\end{equation}

Now note that,
\begin{equation}
\begin{array}{llll}
{\mbf{E}}e^{2\be H_N}{\mbf{1}}_{\{H_N\leq \de N\}} &\leq& \frac{1}{1-4\be^2}& \mbox{ if } \be < \frac{1}{2}\vspace{1ex}\\
&=& \frac{1+\de N}{2}&\mbox{ if }\be=\frac{1}{2}\vspace{1ex}\\
&\leq&\frac{1}{2(2\be-1)}e^{(2\be-1)\de N} & \mbox{ if }
\be>\frac{1}{2}.
\end{array}
\end{equation}

In case $0<\be\leq\frac{1}{2}$, we choose $\de>\log2$ while for
$\frac{1}{2}< \be<1$ we choose $\de$, $\log2< \de
<\frac{\log2}{2\be-1}$ so that by (5) and (6), (4) implies
$$\dsum_{N\geq1}{\mbf{P}}\left[|Z_N'(\be) -
{\mbf{E}}Z_N'(\be)|
> \ep {\mbf{E}}Z_N'(\be)\right] < \infty.$$

Thus, with the choice of $\de$ as specified above, Borel-Cantelli
implies that almost surely eventually,
$$(1-\ep){\mbf{E}}Z_N'(\be)\leq
Z_N'(\be)\leq(1+\ep){\mbf{E}}Z_N'(\be).$$ Combining this with (3) we
have almost surely eventually,
$$(1-\ep){\mbf{E}}e^{2\be H_N}{\mbf{1}}_{\{H_N\leq \de N\}}\leq
\frac{Z_N(\be)}{2^N} \leq(1+\ep){\mbf{E}}e^{2\be
H_N}{\mbf{1}}_{\{H_N\leq \de N\}}.$$

Now fix $K\geq 1$ and $\si \in \{-1, 1\}^K$. Let $Y_N = \dsum_{\si'
\succ \si}e^{-\be H_N(\si')}$, where the sum is over all $\si' \in
\{-1, 1\}^N$ that extend $\si$.

Argument similar to above shows that, with the same $\de$, almost
surely eventually,
$$(1-\ep){\mbf{E}}e^{2\be H_N}{\mbf{1}}_{\{H_N\leq \de
N\}}\leq\frac{Y_N}{2^{N-K}}\leq (1+\ep){\mbf{E}}e^{2\be
H_N}{\mbf{1}}_{\{H_N\leq \de N\}}.$$ As a consequence, almost surely
eventually $\rho_N(\si)$, which by definition is
$\frac{Y_N}{Z_N(\be)}$, lies between $\frac{1-\ep}{1+\ep}\,2^{-K}$
and
$\frac{1+\ep}{1-\ep}\,2^{-K}$ completing the proof.\\

\noindent{\bf Remark 3.1} Clearly for the regime $\be<1$ the
argument shows directly that almost surely $\frac{1}{N}\log Z_N(\be)
\ra \log2$.\\

\noindent{\bf Remark 3.2} Returning to Remark 2.1, if we consider
the environment parametrized by $\al$, it is in general difficult to
evaluate the limiting Gibbs' distribution. However for $\al=2$, our
arguments lead to the convergence of Gibbs' distribution to the
uniform probability, in high temperature regime. (See appendix for
details.)\\

\noindent{\bf Remark 3.2} Hidden in the above argument is a variant
of the strong law of large numbers, which will be taken up
elsewhere.\\

To study the Gibbs' distribution for $\be
>1$, since multiplicative constant cancels out, instead of
$H_N(\si)$ we use the random variables $H_N'(\si) = H_N(\si) + a_N$,
where $a_N = (N-1)\log2$. Mimicking the proof of Lemma 1.2.2 of
Talagrand$^{\cite{T1}}$ yields,
\begin{Prop} For $b \in {\mbf{R}}$,
$$\dlim_{N\ra \infty}{\mbf{P}}( \#\{\si: -H_N'(\si) \geq b\} = k) =
e^{-e^{-b}} \frac{e^{-kb}}{k!}.$$

Moreover, if $\exists$ exactly $k$ many $\si^1, \cdots, \si^k$ in
$\Si_N$ such that $-H_N'(\si^i) \geq b$ then for large $N$ these $k$
points are distributed like $\{X_1, \cdots , X_k\}$ where $X_i$'s
are i.i.d. with density $e^{-(t-b)} {\mbf{1}}_{\{t \geq b\}}$.
\end{Prop}

\noindent {\it Proof:} For fixed $b\in{\mbf{R}}$ and $N$ so large
that $b+a_N>0$, define
\begin{equation}
d_N(b)={\mbf{P}}( -H_N'(\si) \geq b)=\frac{1}{2}\int_b^\infty
e^{-(x+a_N)}dx=\frac{1}{2}e^{-(b+a_N)}.
\end{equation}

By definition of $a_N$ clearly, $2^N d_N(b)=e^{-b}$. Since $\#\{\si:
-H_N'(\si) \geq b\}$ is Binomial with parameters $2^N$ and $d_N$,
Poisson approximation of the binomial completes proof of the first
part.

The last part of the proposition follows from the fact that $H_N'$'s
are i.i.d. and by (7), have density proportional to $e^{-(t+a_N)} {\mbf{1}}_{\{t \geq b\}}$.\\

Let $m>0$ and $\Pi$ be the Poisson point process on $(0,\infty)$
with intensity $x^{-m-1}dx$. Then almost surely, $\Pi$ consists of
summable sequences and hence can be arranged in decreasing order
$(\pi(1),\pi(2),\cdots)$. Let $S(\pi)$ denotes the sum $\sum
\pi(i)$. The distribution of $(\frac{\pi(i)}{S(\pi)}:i\geq 1)$ is
denoted by $PD(m,0)$, called Poisson-Dirichlet distribution with
parameter $m$. For more details and a two parameters family see
Pitman and Yor$^{\cite{PY}}$.
\begin{Prop}
Consider a Poisson point process on ${\mbf{R}}$ with intensity
$e^{-x}dx$, and $(c_i)_{i\geq1}$ an enumeration in decreasing order
of these Poisson points. Then the sequence $$v_i = \frac{e^{\be
c_i}}{\sum_j e^{\be c_j}}$$ has  distribution $PD(\frac{1}{\be},0)$.
\end{Prop}

\noindent {\it Proof:} If $(c_i)$ are Poisson points on ${\mbf{R}}$
with intensity $e^{-x}dx$, $u_i=K e^{\be c_i}$, where $K=\be^\be$,
then $(u_i)$ are Poisson points on $(0,\infty)$ with intensity
$x^{-\frac{1}{\be}-1}dx$. Clearly if $(c_i)$ are enumerated in
decreasing order then so are $(u_i)$. Since $v_i = u_i/\sum_j u_j$,
we conclude that
$(v_i)$ follows $PD(\frac{1}{\be},0)$.\\

Let $\mathcal{S}$ be the set of all decreasing non-negative
sequences with sum at most one. With the $l^1$-metric
$d(\tilde{x},\tilde{y})=\sum |x_i-y_i|$, where $\tilde{x}=(x_i)$ and
$\tilde{y}=(y_i)$, $\mathcal{S}$ is a Polish space. Now let
$\tilde{w}=(w_i)_{i\geq1}$, where $(w_i)_{i\leq 2^N}$ is the
non-increasing enumeration of the (random) Gibbs' weights
$\{G_N(\si) :\si \in \Si_N\}$ with $w_i = 0$ for $i> 2^N$. Let
$\mu_N$ be the law of $\tilde{w}$ on $\mathcal{S}$. Let $\mu$ be the
law of $\tilde{v}=(v_i)_{i\geq 1}$ of Proposition 3.2 on
$\mathcal{S}$.
\begin{theorem}
Let $\be >1$. Then $\mu_N \Rightarrow \mu$ on $\mathcal{S}$, that
is, the law of $\tilde{w}$ converges to $PD(\frac{1}{\be},0)$ as
$N\ra \infty$.
\end{theorem}
\noindent {\it Proof:} One has only to adapt the proof of Theorem
1.2.1 in Talagrand$^{\cite{T1}}$. At the suggestion of the referee
we give a brief outline. Fix $f$ uniformly continuous function on
$\mathcal{S}$, bounded by one and $\ep>0$. Suffices to show $|\int
fd\mu_N-\int fd\mu|<\ep$ for all large $N$.

Let $(c_i)$ be as in Proposition 3.2. Recall that $a_N=(N-1)\log2$
and $H_N'(\si)=H_N(\si)+a_N$. For fixed  $b\in {\mbf{R}}$, Put
temporarily,
$$\begin{array}{lll}
Z_N(\be)&=&\sum_\si e^{-\be H_N'(\si)},\vspace{1ex}\\
Z_N(\be,b) &=& \sum_\si e^{-\be H_N'(\si)}{\mbf{1}}_{\{-H_N'(\si)
\geq b\}},\vspace{1ex}\\
w_i^b &=& \frac{e^{\be h_i}}{Z_N(\be, b)}{\mbf{1}}_{\{h_i\geq b\}}
\mbox{ if } w_i=\frac{e^{\be h_i}}{Z_N(\be)},\vspace{1ex}\\
v_i^b&=&\frac{e^{\be c_i}{\mbf{1}}_{\{c_i\geq b\}}}{\sum e^{\be
c_i}{\mbf{1}}_{\{c_i\geq b\}}}.
\end{array}$$
Denote $\tilde{w}^b=(w_i^b)_{i\geq1}$ and
$\tilde{v}^b=(v_i^b)_{i\geq1}$. Let $\mu_N^b$ (respectively,
$\mu^b$) be the law of $\tilde{w}^b$ (respectively, $\tilde{v}^b$)
on $\mathcal{S}$. Observe,
\begin{equation}
\begin{array}{lll}
|\int fd\mu_N-\int fd\mu|&<&|\int fd\mu_N-\int fd\mu_N^b|\,+\,|\int
fd\mu_N^b-\int fd\mu^b|\vspace{1ex}\\
&& \hspace{23ex}+|\int fd\mu^b-\int fd\mu|.
\end{array}
\end{equation}

{\bf Firstly}, we show that for given $\de>0$, there exists $b_0$
(depending on $\be$ and $\de$) such that for $b\leq b_0$,
\begin{equation}
\dsty\limsup_{N\ra \infty} {\mbf{P}}\left(\frac{Z_N(\be)
-Z_N(\be,b)}{Z_N(\be)} \geq \de\right) \leq \de.
\end{equation}

Then, using $d(\tilde{w},{\tilde{w}^b})=2\frac{Z_N(\be)-Z_N(\be,
b)}{Z_N(\be)}$ and that $f$ is bounded and uniformly continuous, the
first term on the right side of (8) can be made $<\frac{\ep}{3}$ for
all large $N$. To see (9), we proceed as follows:

Since $Z_N(\be) \leq e^{\be x} \Ra \{ \#\{\si : - H_N(\si) \geq x\}
=0\}$, by Proposition 3.1, $\dlim_{N\ra \infty} {\mbf{P}}( \#\{\si:
-H_N(\si) \geq
 x\}=0)=e^{-e^{-x}},$ so there exists $\eta >0$ such that  for large $N$,
\begin{equation}
{\mbf{P}}(Z_N(\be) \leq \eta) \leq \frac{\de}{2}.
\end{equation}

Again, for fixed $x\in {\mbf{R}}$ with $N$ so large that $x+a_N>0$,
${\mbf{E}}(Z_N(\be) -Z_N(\be,x))= \frac{1}{\be-1}e^{(\be-1)x} -
\frac{2}{\be^2 - 1} e^{-(\be-1) a_N}.$ Since $\be>1$ this can be
made sufficiently small by an appropriate choice of large negative
quantity $x$ and large $N$. So that, by Chebycheff's inequality, we
can get $b_0$ (depending on $\eta, \de, \be$) such that for $b\leq
b_0$ and large $N$,
\begin{equation}
{\mbf{P}}[Z_N(\be) -Z_N(\be,b)\geq \eta\de] \leq  \frac{\de}{2}.
\end{equation}

Now (10) and (11) imply (9).

{\bf Secondly}, the last term in (8) can be made small by choosing
$b$ large negative quantity since $\mu^b \Ra \mu$ as $b \ra
-\infty$. Fix now such a number $b$. There is no loss to assume (9)
also holds.

{\bf Finally}, the middle term in (8) can be made arbitrary small by
choosing $N$ large, since by the last part of Proposition 3.1
$\mu_N^b \Ra \mu^b$.

Thus for large $N$ each of the three terms on the right hand side of
(8) can be made smaller than $\frac{\ep}{3}$.

\section*{Appendix}
\noindent{\bf Proposition:} {\em Consider  REM with $H_N\sim$
Gaussian $(0,N)$ and fix $0<\be < \sqrt{2\log2}$. Then almost
surely, for any $K \geq 1$ and any $\si \in \{-1, 1\}^K$,
$\rho_N(\si) \ra \frac{1}{2^K}$ as $N \ra \infty$, where $\rho_N$ is
the marginal of $G_N$ on $\{-1, 1\}^K$.}

A stronger version with a difficult proof is in
Talagrand$^{\cite{T1}}$. The purpose of this appendix is to explain
how the techniques used in Theorem 3.1 above apply to this case.

For $0<\be<\sqrt{2\log2}$, define $Z_N'(\be)=\sum_\si e^{-\be
H_N(\si)}{\mbf{1}}_{\{-H_N(\si)\leq \de N\}}$, where $\de$ will be
chosen latter depending on $\be$.

Since ${\mbf{P}}\{-H_N(\si)> \de N\} \leq \frac{1}{\sqrt{2\pi
N}}e^{-\frac{\de^2N}{2}}$, for any $\de>\sqrt{2\log2}$,
Borel-Cantelli implies that almost surely, eventually
\begin{equation}
Z_N(\be)=Z_N'(\be).
\end{equation}

Just as in Theorem 3.1, we have,
\begin{equation}{\mbf{P}}\left[|Z_N'(\be) -
{\mbf{E}}Z_N'(\be)|
> \ep {\mbf{E}}Z_N'(\be)\right] < \frac{{\mbf{E}}e^{2\be
H_N}{\mbf{1}}_{\{H_N\leq \de N\}}}{\ep^2 2^N ({\mbf{E}}e^{\be
H_N}{\mbf{1}}_{\{H_N\leq \de N\}})^2}.
\end{equation}

Since $\de>\sqrt{2\log2}>\be$, \begin{equation}{\mbf{E}}e^{\be
H_N}{\mbf{1}}_{\{H_N\leq \de N\}}=\frac{e^{\frac{\be^2
N}{2}}}{\sqrt{2\pi N}}\int_{-\infty}^{\de N}e^{\frac{(x-\be N)^2
}{2N}}dx>\frac{1}{2}e^{\frac{\be^2 N}{2}}.\end{equation}

Now note that, $${\mbf{E}}e^{2\be H_N}{\mbf{1}}_{\{H_N\leq \de
N\}}=\frac{e^{2\be^2 N}}{\sqrt{2\pi N}}\int_{(2\be-\de)N}^\infty
e^{-\frac{x^2}{2N}}dx.$$ Thus,
\begin{equation}
\begin{array}{llll}
{\mbf{E}}e^{2\be H_N}{\mbf{1}}_{\{H_N\leq \de N\}} &\leq& e^{2\be^2
N}& \mbox{ if } \be \leq \frac{\de}{2}\\
&\leq& \frac{1}{(2\be-\de)\sqrt{2\pi N}}
e^{(2\de\be-\frac{\de^2}{2})N}&\mbox{ if }\be>\frac{\de}{2}.
\end{array}
\end{equation}

In case $0<\be<\sqrt{\log2}$, we choose $\de=2\sqrt{\log2}$ so that
$\be\leq\frac{\de}{2}$ and hence by (14) and (15), (13) implies
$$\dsum_{N\geq1}{\mbf{P}}\left[|Z_N'(\be) -
{\mbf{E}}Z_N'(\be)|
> \ep {\mbf{E}}Z_N'(\be)\right] < \frac{4}{\ep^2}\dsum_{N\geq1}e^{-N(\log2 - \be^2)}<\infty.$$

In case $\sqrt{\log2}\leq\be < \sqrt{2\log2}$, one observes that
$\sqrt{2\log2}<2\be -\sqrt{2(\be^2-\log2)}$ so that we can choose
$\de \in (\sqrt{2\log2},2\be -\sqrt{2(\be^2-\log2)})$. With such a
choice, $\be>\frac{\de}{2}$ and again by (14) and (15), (13) implies
$$\dsum_{N\geq1}{\mbf{P}}\left[|Z_N'(\be) -
{\mbf{E}}Z_N'(\be)|
> \ep {\mbf{E}}Z_N'(\be)\right] <\infty.$$
Thus in either case, if $\de$ is chosen as specified, Borel-Cantelli
implies that almost surely eventually,
$$(1-\ep){\mbf{E}}Z_N'(\be)\leq
Z_N'(\be)\leq(1+\ep){\mbf{E}}Z_N'(\be)$$ and the proof is completed
by repeating the same argument of Theorem 3.1.
\section*{Acknowledgments}

The author would like to thank B. V. Rao for valuable discussions.
Also we would like to thank the referees for their helpful comments.



\begin{thebibliography}{99}
\footnotesize
\bibitem{B}
{\sc Bovier, A.} (2001). Statistical Mechanics of Disordered
Systems, {\em MaPhySto Lecture Notes}, 192pp.

\bibitem{CEGG}
{\sc Contucci, P., Esposti, M.~D., Giardin\`{a}, C. and Graffi, S.}
(2003). Thermodynamical Limit for Corelated Gaussian Random Energy
Models. {\em Commun. Math. Phys.}, {\bf 236}, 55--63.

\bibitem{DZ}
{\sc Dembo, A. and Zeitouni, O.} (1998). {\em LARGE DEVIATION:
Techniques and Applications}, Second Edition, Springer-Verlag, New
York.

\bibitem{D}
{\sc Derrida, B.} (1981). Random Energy Model: An Exactly Solvable
Model of Disordered Systems. {\em Phys. Rev.}, {\bf B24},
2613--2626.

\bibitem{DW}
{\sc Dorlas, T.~C. and Wedagedera, J.~R.} (2001). Large Deviations
and The Random Energy Model. {\em Int. J. of Mod. Phy. B}, {\bf 15,}
No. 1, 1--15.

\bibitem{OP}
{\sc Olivieri, E. and Picco, P.} (1984). On existance of
Thermodynamics for the Random Energy Model. {\em Commun. Math.
Phys.}, {\bf 96}, 125--144.

\bibitem{PY}
{\sc Pitman, J. and Yor, M.} (1997). The two-parameter
Poissin-Dirichlet distribution derived from a stable subordinator.
{\em Ann. Probab.}, {\bf 25}, 855--900.

\bibitem{T1}
{\sc Talagrand, M.} (2003). {\em SPIN GLASSES: A Challenge for
Mathematicians}, Spinger-Verlag, New York.

\end{thebibliography}
\end{document}